# Optimal Distributed Generation Planning in Active Distribution Networks Considering Integration of Energy Storage

Yang Li [a], Bo Feng [a], Guoqing Li [a], Junjian Qi [b], Dongbo Zhao [c], Yunfei Mu [d]

[a] School of Electrical Engineering, Northeast Electric Power University, Jilin 132012, P.R. China
[b] Department of Electrical and Computer Engineering, University of Central Florida, Orlando, FL 32816 USA
[c] Energy Systems Division, Argonne National Laboratory, Argonne, IL 60439 USA
[d] Key Laboratory of Smart Grid of Ministry of Education, Tianjin University, Tianjin 300072, China

*Abstract*—A two-stage optimization method is proposed for optimal distributed generation (DG) planning considering the integration of energy storage in this paper. The first stage determines the installation locations and the initial capacity of DGs using the well-known loss sensitivity factor (LSF) approach, and the second stage identifies the optimal installation capacities of DGs to maximize the investment benefits and system voltage stability and to minimize line losses. In the second stage, the multi-objective ant lion optimizer (MOALO) is first applied to obtain the Pareto-optimal solutions, and then the 'best' compromise solution is identified by calculating the priority memberships of each solution via grey relation projection (GRP) method, while finally, in order to address the uncertain outputs of DGs, energy storage devices are installed whose maximum outputs are determined with the use of chance-constrained programming. The test results on the PG&E 69-bus distribution system demonstrate that the proposed method is superior to other current state-of-the-art approaches, and that the integration of energy storage makes the DGs operate at their pre-designed rated capacities with the probability of at least 60% which is novel.

*Index Terms*--active distribution network; distributed generation planning; two-stage optimization; energy storage; chance-constrained programming; multi-objective ant lion optimizer.

**Nomenclature**

| | |
|---|---|
| $\alpha_i$ | loss sensitivity factor of real power loss at the $i$th bus |
| $P_L$ | total system active power losses; |
| $Q_L$ | total system reactive power losses; |
| $U_i \angle \delta_i$ | complex voltage at the bus $i$; |
| $R_{ij}$ | resistance between the buses $i$ and $j$; |
| $P_i$ | active power injections at the bus $i$ |
| $Q_j$ | reactive power injections at the bus $j$ |
| $N$ | number of buses |
| $C_i^{\text{GP}}$ | on-grid price of DGs at the bus $i$ |
| $C_i^{\text{GS}}$ | government subsidy of DGs at the bus $i$ |
| $S_i^{\text{rated}}$ | rated capacity of DGs at the bus $i$ |
| $\lambda_i^{\text{CF}}$ | capacity factor of DGs at the bus $i$ |
| $C_i^{\text{MC}}$ | maintenance cost of DGs at the bus $i$ |
| $C_i^{\text{FIC}}$ | fixed investment cost of DGs at the bus $i$ |
| $\xi_i^{\text{DG}}$ | annual conversion factor of fixed investment cost of DGs at the bus $i$ |
| $C_P$ | investment benefit of DGs per year |
| $C_T$ | annual revenue of investment in DGs |
| $C_I$ | annual cost of investment in DGs |
| $VSF_{i+1}$ | voltage stability factor of the $i+1$ bus |
| $VSF_{total}$ | voltage stability factor of the total distribution network |
| $R_m$ | equivalent resistance of the branch $m$ |
| $P_{\text{loss}}(m)$ | active power loss of the branch $m$ |
| $P_{\text{DG},i+1}$ | active power output of DGs at the bus $i+1$ |
| $Q_{\text{DG},i+1}$ | reactive power output of DGs at the bus $i+1$ |
| $P_{\text{loss,total}}$ | total line losses of the distribution system |
| $NB$ | number of branches in the network |
| $U_{\min}$ | lower bounds of the voltage at the bus $i$ |

| | | | |
|---|---|---|---|
| $U_{max}$ | upper bounds of the voltage at the bus $i$ | $v_{out}$ | cut-out wind speed |
| $P_{Di}$ | active power of load at the bus $i$ | $\eta$ | conversion efficiency of PV units |
| $Q_{Di}$ | reactive power of load at the bus $i$ | $r_{max}$ | maximum solar-irradiance intensity |
| $P_{LO}$ | system active power losses without DGs | *Abbreviation* | |
| $Q_{LO}$ | system reactive power losses without DGs | ADNs | active distribution networks |
| $S_{Lm}$ | actual complex power in the branch $m$ | DG | distributed generation |
| $S_{Lm(rated)}$ | rated complex power in the branch $m$ | MINLP | mixed integer non-linear programming |
| $P_{Swing}$ | active power of the swing bus | CCP | chance-constrained programming |
| $Q_{Swing}$ | reactive power of swing bus | WT | wind turbine |
| $NI$ | number of indications | PV | photovoltaic |
| $\gamma_{lk}$ | grey relation coefficient between the $k$th indication in the $l$th scheme | MT | micro-gas turbine |
| | | ALO | ant lion optimizer |
| $w_k$ | weight of each indication in the scheme | MOEAs | multi-objective evolutionary algorithms |
| $D_l$ | priority membership of the $l$th scheme | MOALO | multi-objective ant lion optimizer |
| $P_w$ | optimal output of wind power energy | NSGA | non-dominant sorting genetic algorithm |
| $P_s$ | optimal output of solar energy | MOPSO | multi-objective particle swarm optimization |
| $P_t^W$ | output of wind power at the $t$th sampling | | |
| $P_t^S$ | output of solar energy at the $t$th sampling | MOHS | multi-objective harmony search |
| $\omega$ | given confidence level | LSF | loss sensitivity factor |
| $v_r$ | rated wind speed | ESDs | energy storage devices |
| $v_{in}$ | cut-in wind speed | PPF | probabilistic power flow |

## I. INTRODUCTION

Nowadays, with the increasingly high penetration of renewable distributed generation (DG) sources, active distribution networks (ADNs) have been regarded as an important solution to achieve power system sustainability and energy supply security [1, 2]. Recently, it is becoming an inevitable trend to make full use of renewable DGs such as wind turbines (WT) and photovoltaic (PV) units, since they have substantial advantages, such as power loss reduction, greenhouse gas emission reduction, flexible voltage regulation, peak-load shaving, higher power quality, supply reliability enhancement [3, 4], etc. However, extensive penetration of DGs greatly increases the risks of safe and economic operation of distribution networks since renewable DGs have inherently intermittent nature [5, 6], which makes the planning more challenging than ever before. The traditional distribution network planning options, such as the addition or expansion of substations and lines, are unable to meet the needs of modern complex ADNs facing all alternatives together with generation and load uncertainties [6, 8, 9]. Therefore, it is necessary to deal with such key challenges in the issue of optimal DG placement.

For this problem, many attempts have been made in the past two decades to solve it using different methodologies, including analytical approach [10-13], mixed integer non-linear programming (MINLP) [14-16], Kalman filter algorithm [17], and computational intelligence approach [18-20]. An analytical method was firstly introduced in [10] to find out the optimal placement of a single DG in both radial and meshed networks to minimize power losses. However, this approach only optimizes siting and considers DG sizing as fixed. In [11], the locating and sizing of DGs are identified instantaneously using an analytical strategy. In [12], an improved analytical approach for multiple DG planning was proposed for reducing energy losses, and the method was examined on three different test systems. A simple analytical type approach was recently presented to optimize the loss related to the active and reactive components of DG branch current. Different from the aforesaid works, a methodology based on based on MINLP was developed in [14] for optimally planning different types of renewable DGs for minimizing annual power loss. In [15], same approach has been implemented for optimal renewable DG placement and sizing to improve the voltage stability margin. Considering the uncertainties of loads and DG outputs, a technique based on multi-objective MINLP was proposed for benefit maximization in distribution systems in [16]. In [17], the optimal DG placement is determined via a Kalman filter algorithm to minimize the losses. More recently, computational intelligence techniques have been successfully adopted to deal with the DG planning issues. In [18], a genetic algorithms-based expansion planning model considering DG integration and conventional alternatives for expansion was presented for ADNs, and two multiple scenarios analysis approaches were employed to tackle uncertainties of DG and load response. In [19], a method for optimum DG siting and sizing of multi-DG units based on maximization of system loadability is proposed with the use of hybrid particle swarm optimization. The work in [20] addressed a multi-objective formulation for simultaneous allocation of distributed energy resources to maximize annual savings by using improved particle swarm optimization.

To sum up, DG planning methodologies for distribution networks have been widely studied, but there is relatively little focus on considering the generation and load uncertainties [6]. Furthermore, due to the high penetration of renewable DG resources, the uncertainties in ADNs are becoming significantly larger than traditional networks [21, 22], which causes the actual power outputs of DGs to be barely achieving their pre-designed rated capacities by using the aforesaid methods. In general, it allows obtaining a planning scheme more committed with practical operations to take into account uncertainties during the optimization process in the planning stage, since typical scenarios with their occurrence probabilities are analyzed [18].

Recent research findings have shown that energy storage plays an increasingly important role in optimal DG allocation in distribution networks for the purpose of integrating intermittent renewable generation and loads [21-24], since energy storage devices (ESDs) can effectively shift energy generation and consumption across time spots [25]. After years of research and practice, there are a large set of storage technologies available to support renewable DGs [26], such as battery energy storage [27], supercapacitors [28-31], fuel cells [26, 32], etc. In addition, rapid advances recently made in the area of energy storage [32-35] provide more powerful and flexible supports for renewable DG integration.

The contributions of this paper are threefold as follows. (1) A two-stage optimization approach for DG planning in distribution networks is proposed. (2) More importantly, different from most existing works on DG planning in which the power fluctuations from DGs are ignored, we consider the random power output characteristics of DGs and determine the maximum outputs of energy storage devices (ESDs) through the utilization of chance-constrained programming (CCP) to make sure that the DG power outputs can achieve their pre-designed rated capacities with high probability under real operation condition. (3) Last but not least, the impacts of energy storage integration are deeply analyzed with the use of probabilistic power flow (PPF) calculations.

The remainder of this paper is structured as follows. Section II gives the problem formulation. A detailed description of the two-stage model-solving process is put forward in Section III. Section IV examine effectiveness of the proposal on the PG&E 69-bus distribution system. Finally, the conclusions are drawn from the simulation results in Section V.

## II. PROBLEM FORMULATION

In this section, we define the objective functions and technical constraints of the multi-objective optimization model for determining the optimal locations and sizes of DGs.

### A. The objective functions

#### 1) Investment benefit

As is known, investment in DG units is a very attractive distribution planning option [36, 37]. In this work, the DG investment benefit is defined as the ratio of the annual income and the investment of DGs.

$$C_p = C_T / C_I, \tag{1}$$

where $C_T$ is the annual revenue of investment in DGs, including the benefits of selling electricity and the policy subsidy of DG; $C_I$ is the annual cost of investment in DGs, including installation costs, operation and maintenance costs, and fuel costs of DG.

$$C_T = 8760 \sum_{i=1}^{N} (C_i^{GP} + C_i^{GS}) S_i^{rated} \lambda_i^{CF}$$

$$C_I = 8760 \sum_{i=1}^{N} C_i^{MC} S_i^{rated} \lambda_i^{CF} + \sum_{i=1}^{N} C_i^{FIC} S_i^{rated} \xi_i^{DG}$$

where $N$ is the number of total buses, $C_i^{GP}$ and $C_i^{GS}$ are respectively the on-grid price and the government subsidy of DGs at the bus $i$; $S_i^{rated}$ is the rated capacity of DGs at the bus $i$; $\lambda_i^{CF}$ is the capacity factor of DGs at the bus $i$; $C_i^{MC}$ and $C_i^{FIC}$ are respectively the maintenance cost and fixed investment cost of DGs at the bus $i$; $\xi_i^{DG}$ is the annual conversion factor of fixed investment cost of DGs at the bus $i$.

*2) Voltage stability factor*

An efficient and simple voltage stability factor $VSF_{total}$ is recently derived to determine voltage stability condition of the whole distribution network in [10]. Based on the power flow calculation, the voltage stability factor of the $(i+1)$ th bus is obtained as [13]:

$$VSF_{i+1} = 2U_{i+1} - U_i. \tag{2}$$

By summing the voltage stability factors of all the load buses, the voltage stability of the whole system can be justified by the indictor $VSF_{total}$. The lower value of $VSF_{total}$ suggests more voltage unstable operation. The indictor $VSF_{total}$ is defined as follows.

$$VSF_{total} = \sum_{i=1}^{N-1}(2U_{i+1} - U_i) \tag{3}$$

*3) Line loss*

For a given branch m connected between the buses $i$ and $i+1$, as shown in Fig. 1.

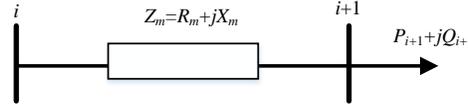

**Fig. 1.** Branch model in radial distribution systems

The active power loss in branch *m* is specified by $I_m^2 R_m$, which can be given by

$$P_{loss}(m) = R_m \left[\frac{(P_{i+1} - P_{DG,i+1})^2 + (Q_{i+1} - Q_{DG,i+1})^2}{U_{i+1}^2}\right] \tag{4}$$

where $R_m$ is the equivalent resistance of branch *m* connecting bus $i$ and bus $i+1$; $P_{DG,i+1}$ is the active power output of DGs at bus $i+1$; $Q_{DG,i+1}$ is the reactive power output of DGs at bus $i+1$; $U_{i+1}$ is the voltage magnitude of bus $i+1$.

Accordingly, the line losses of the whole distribution system with integration of DGs is expressed as

$$P_{loss\_total} = \sum_{m=1}^{NB} P_{DG\_loss}(m) \tag{5}$$

where *NB* is the number of branches in the network.

B. *Technical constraints*

The voltage at each bus is limited by

$$U_{min} \leq |U_i| \leq U_{max} \tag{6}$$

where $U_{min}$ and $U_{max}$ are, respectively, the lower and upper bounds of the voltage at the bus *i*.

The maximum penetration of DGs defined as the total capacity of DG units is limited as

$$\sum_{i=2}^{N} P_{DG,i} \leq \sum_{i=2}^{N} P_{Di} + P_{LO} \tag{7}$$

$$\sum_{i=2}^{N} Q_{DG,i} \leq \sum_{i=2}^{N} Q_{Di} + Q_{LO} \tag{8}$$

where $P_{Di}$ and $Q_{Di}$ are respectively active and reactive power of load at the bus *i*; and $P_{LO}$ and $Q_{LO}$ are total system active and reactive power losses without DG units, respectively.

The complex power through any line must be less than its rating value as:

$$S_{Lm} \leq S_{Lm(rated)} \tag{9}$$

where $S_{Lm}$ and $S_{Lm(rated)}$ are the actual and the rated complex power in branch *m*.

The total power consumption should be equal to the total power supply at each bus:

$$P_{Swing} + \sum_{i=2}^{N} P_{DG,i} = \sum_{i=2}^{N} P_{Di} + P_L \tag{10}$$

$$Q_{Swing} + \sum_{i=2}^{N} Q_{DG,i} = \sum_{i=2}^{N} Q_{Di} + Q_L \tag{11}$$

where $P_{Swing}$ and $Q_{Swing}$ are, respectively, active and reactive power of the swing bus, and $Q_L$ is the total system reactive power loss.

### III. MODEL SOLVING

The first stage determines the locations and initial capacities of DGs according to the loss sensitivity factor (LSF) approach; the second stage determines the optimal installation capacities of DGs using MOALO and GRP. Finally, to deal with the randomness of the outputs of DGs, energy storage devices are installed whose maximum outputs are determined by using CCP.

#### A. 1st stage: loss sensitivity factor

In the first stage, the location and initial capacity of DGs are determined through the use of the well-known loss sensitivity factor (LSF) approach, which has been widely used for solving the problems of optimal capacitor and DG allocations [10, 15]. The total real power loss $P_L$ in a power system is represented by [15]

$$P_L = \sum_{i=1}^{N}\sum_{j=1}^{N}[\alpha_{ij}(P_iP_j + Q_iQ_j) + \beta_{ij}(Q_iP_j - P_iQ_j)] \quad (12)$$

where

$$\alpha_{ij} = \frac{R_{ij}}{U_iU_j}\cos(\delta_i - \delta_j), \quad \beta_{ij} = \frac{R_{ij}}{U_iU_j}\sin(\delta_i - \delta_j)$$

Here, $P_i$ and $P_j$ are respectively the active power injection at the buses $i$ and $j$, $Q_i$ and $Q_j$ are the reactive power injections at the buses $i$ and $j$, $U_i \angle \delta_i$ denotes complex voltage at the bus $i$, $R_{ij}$ is resistance between the buses $i$ and $j$.

And thereby, the LSF of real power loss at the $i$th bus $\alpha_i$ is:

$$\alpha_i = \frac{\partial P_L}{\partial P_i} = 2\sum_{j=1}^{N}(\alpha_{ij}P_j - \beta_{ij}Q_j) \quad (13)$$

The specific procedures to find the optimal DG locations and initial capacities of DG units are as follows [15]:

Step 1: Input the total number of DG units to be installed.

Step 2: Calculate losses according to (12) by solving the load flow for the base case.

Step 3: Find the optimal DG location. First, calculate the LSF at each bus according to (13), then choose the bus with the largest LSF value as the highest priority bus.

Step 4: Find the optimal size of DG and calculate losses with the minimum loss by running load flow.

Step 5: Update load data to allocate the next DG.

Step 6: Terminate if the constraint is unsatisfied; otherwise, repeat steps 2 to 5.

#### B. 2nd stage: multi-objective optimization and decision-making

1) *Multi-objective ant lion optimizer*

Ant Lion Optimizer (ALO) proposed by Mirjalili in 2015 is a novel nature-inspired algorithm [38], and it has been widely used for solving various engineering problems [39, 40]. The ALO mimics the hunting mechanism of ant lions in nature. An ant lion larva digs a cone shaped pit in the sand by moving along a circular path and throwing out sands with its massive jaw. After digging the trap, the larva hides underneath the bottom of the cone and waits for insects to be trapped in the pit. The edge of the cone is sharp enough for insects to fall to the bottom of the trap easily. Once the ant lion realizes that a prey is in the trap, it tries to catch it. Then, the prey is pulled under the soil and consumed. After consuming the prey, ant lions throw the leftovers outside the pit and prepare the pit for the next hunt [38]. Five main steps of hunting prey, the random walk of ants, building traps, entrapment of ants in traps, catching preys, and re-building traps are implemented.

In order to solve multi-objective optimization problems, Multi-objective ALO (MOALO) was designed with equipping ALO with an archive and ant lion selection mechanism based on Pareto optimal dominance [41]. Large amounts of benchmarking tests demonstrate that MOALO are superior to other currently popular multi-objective evolutionary algorithms (MOEAs) such as multi-objective particle swarm optimization (MOPSO) and Non-dominated Sorting Genetic Algorithm-II (NSGA-II), no matter whether solution quality or convergence speed [41].

The optimization process of the MOALO algorithm mainly includes the following steps:

Step 1: Initialize all ants and antlions, and set the current iteration number *gen* and the ant number *num_ant* are equal to 1;

Step 2: Select a random antlion from the archive;

Step 3: Select the elite using Roulette wheel from the archive;

Step 4: Create a random walk and normalize it;

Step 5: Update the position of ant;

Step 6: If every ant has been traversaled, then go to step 7; otherwise, go to step 2, and the *num_ant* is incremented by 1;

Step 7: Calculate the objective values of all ants;

Step 8: Update the archive;

Step 9: Judge whether the termination condition is met or not. If the current iteration exceed the pre-defined maximum of iterations, then go to Step 10; otherwise, go to step 2, and the *gen* is incremented by 1.

Step 10: output the Pareto-optimal solutions.

2) *Grey Relation Projection*

After having obtained the set of Pareto optimals, it is an important work to determine the best compromise solution in decision making process. GRP theory is a powerful tool for analyzing the relationship between sequences with grey information and has been successfully applied in a variety of fields [42].

According to the theory, the projection $V_l$ of the *l*th scheme on the ideal reference scheme is given by

$$V_l^{+(-)} = \sum_{k=1}^{NI} \gamma_{lk}^{+(-)} \frac{w_k^2}{\sqrt{\sum_{k=1}^{NI}(w_k)^2}} \qquad (14)$$

where superscript "+"/"-" respectively denote positive/negative schemes, $NI$ is the number of indications, $\gamma_{lk}$ is the grey relation coefficient between the *k*th indication in the *l*th scheme, $w_k$ is the weight of each indication in the scheme. Here, for the convenience of description, the weights of all objective functions are assumed to be equal, but they can be adjusted according to decision-makers' preferences. The priority membership $D$ of the *l*th scheme $D_l$ is

$$D_l = \frac{(V_0 - V_l^-)^2}{(V_0 - V_l^-)^2 + (V_0 - V_l^+)^2}, \quad 0 \leq D_l \leq 1 \qquad (15)$$

where $V_0$ is equal to $V_l$ when $\gamma_{lk}$ is taken as 1. The equation shows that the higher the priority membership is, the better the scheme will be. Therefore, the results with the highest priority membership will be chosen as the 'best' compromise solutions.

3) *Determination of energy storage maximum output*

Traditional distribution network planning is usually based on the ideal system operating conditions, without considering the uncertainty of DG outputs and loads [5, 6]. However, due to the random output characteristics of WT and PV units, the actual power outputs of DGs are hardly achieving (with a very low or even zero probability) their pre-designed rated capacities in the planning stage by using traditional planning methodologies. In addition, compared with the intermittent power outputs of DGs, the fluctuations of loads are relatively much smaller [22]. Taking into consideration the above factors, we have simplified this problem by only considering the uncertainty of DGs (ignoring the load uncertainties) in this work.

Recent work has shown that spinning reserve is an important resource to reduce operation costs and maintain the system's reliability by means of compensating the power fluctuations of renewable DGs [43]. With the rapid development of energy storage technology, it has become an effective means for meeting spinning reserve requirements [26, 32]. Taking into account that batteries have advantages of easy scaling in terms of power rating and capacity [27], for the sake of simplicity, a general type of battery is utilized for providing spinning reserve services in this work. The ESDs should be able to compensate the difference between the optimal output and the stochastic output of DGs. Given a pre-assigned confidence level, the output of ESDs can be calculated by using the CCP theory. Here, the used CCP model is

$$\begin{aligned} &\min \ P_{\text{Reest}} \\ &P\{P_{\text{Reest}} \geq P_w + P_s - P_t^w - P_t^s\} \geq \omega, \ \forall t \end{aligned} \qquad (16)$$

where $R_e = P_w + P_s - P_t^w - P_t^s$, $P_{Reest}$ is the energy storage maximum output to be determined, $P\{\cdot\}$ is the probability of an event, $P_w$ and $P_s$ are the optimal output of the WT and PV unit, $P_t^W$ is the output of the WT unit at the $t$th sampling, $P_t^S$ is the output of the PV unit at the $t$th sampling, and $\omega$ is the given confidence level.

For a given confidence level $\omega$, the probability distribution of $P_t^W$ and $P_t^S$ are calculated by using Monte Carlo simulation techniques, $P_w$ and $P_s$ are extracted from the best comprise solution, and then the energy storage maximum output $P_{Reest}$ can be determined through the use of the CCP theory.

*C. Solving process*

The proposed problem-solving process consists of the following two stages: the first stage determines the installation locations and the initial capacities of DGs via the LSF approach, and the second stage identifies the optimal installation capacities of DGs through the use of MOALO and GRP.

The specific solving process based on the two-stage optimization is shown in Fig. 2.

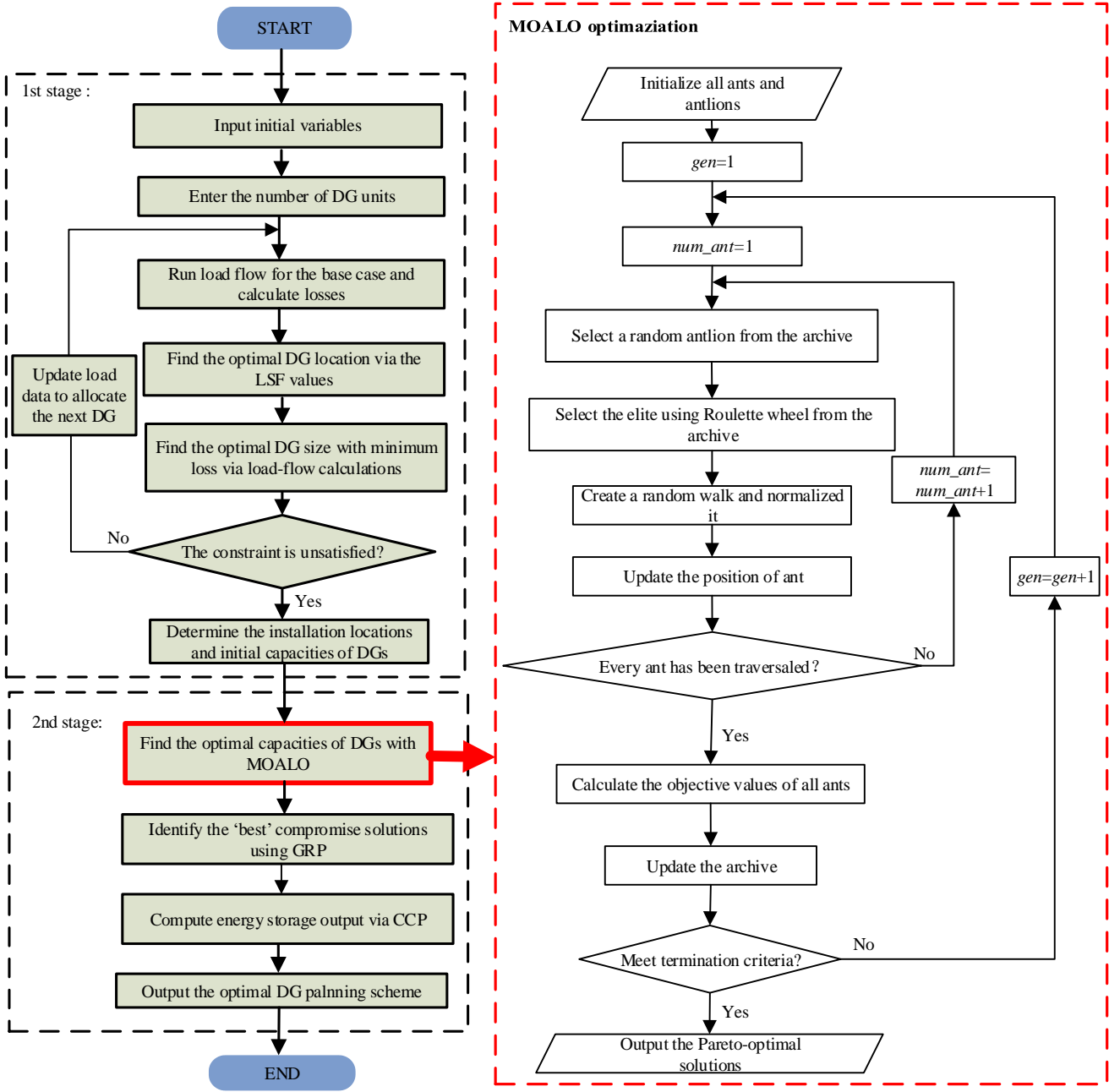

**Fig. 2.** Flowchart of the proposed approach

## IV. CASE STUDIES

### A. PG&E69-bus system

The PG&E 69-bus distribution system is a well-known test case for distribution network planning [8, 9, 11, 12, 19, 20, 37, 39], as shown in Fig. 3. In this system, the voltage level is 12.66 kV, the total active and reactive loads are, respectively, 3715 kW and 2300 kVAr. This network is connected to external power grid via the bus 1 (swing bus) and the system parameters can be found in [19].

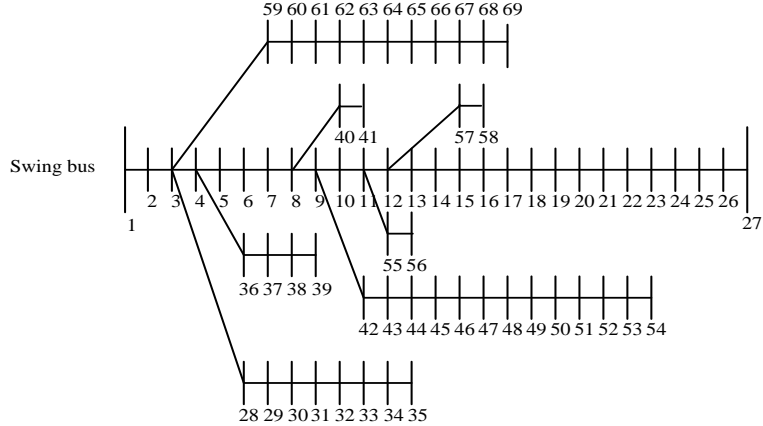

**Fig. 3.** Network diagram of PG&E 69-bus system

### B. Parameter settings

The operating parameters of the wind power and photovoltaic generation are listed in Table I. Where $v_r$ denotes the rated wind speed, $v_{in}$ represents the cut-in wind speed, $v_{out}$ is the cut-out wind speed, $\eta$ is the conversion efficiency of PV units, $r_{max}$ is the maximum solar-irradiance intensity.

TABLE I. PARAMETERS OF WIND POWER AND PHOTOVOLTAIC GENERATION

| $v_{in}$ (m/s) | $v_r$ (m/s) | $v_{out}$ (m/s) | $\eta$ | $r_{max}$ (W/m$^2$) |
|---|---|---|---|---|
| 4 | 16 | 28 | 12 | 28 |

The parameters of investment are listed in Table II. Here, WT, PV, and MT are, respectively, the wind power unit, the photovoltaic unit and the micro-gas turbine. Evidently, the DG costs and profits per kWh are different in different countries or regions, and they may vary due to technology progress and policy changes.

TABLE II. INVESTMENT PARAMETERS OF DGS

| DG | Investment cost $C^{FIC}$ (×10$^4$ $/kWh) | Maintenance cost $C^{MC}$ ($/kWh) | On-grid price $C^{GP}$ ($/kWh) | Government subsidy $C^{GS}$ ($/kWh) | Conversion factor $\gamma^{DG}$ | Capacity factor $\lambda^{CF}$ |
|---|---|---|---|---|---|---|
| WT | 0.163 | 0.0047 | 0.08 | 0.036 | 0.1006 | 0.35 |
| PV | 0.667 | 0.0019 | 0.08 | 0.036 | 0.0843 | 0.29 |
| MT | 0.164 | 0.0283 | 0.064 | 0 | 0.1006 | 1.00 |

### C. Optimal DG allocation calculation

The locations and the initial capacities of DGs obtained from the first stage are listed in Table III. The results are obtained based on the LSF values at different buses. Note that the location of DGs will not change in subsequent calculations while the capacity will be updated in the iterations to achieve optimal results.

TABLE III. INSTALLATION POSITIONS AND CAPACITIES OF DG

| DG installation locations | Capacities of corresponding DGs /kW |
|---|---|
| bus 64 | 47.7 |
| bus 49 | 114.5 |
| bus 50 | 170.9 |
| bus 61 | 1533.5 |

The Pareto optimal solution set of optimal allocation of DGs is shown in Fig. 4, and the best comprise solution (marked as solution A) and a set of typical solutions are shown in Table IV.

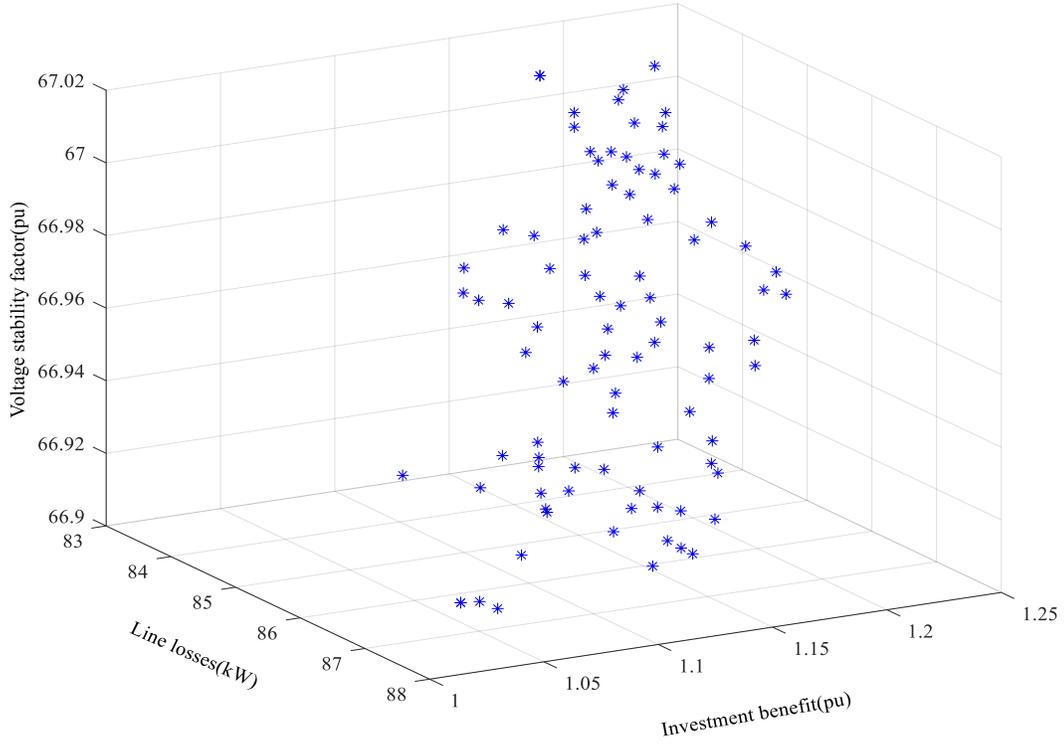

**Fig. 4.** Distribution of Pareto optimal solutions

TABLE IV. SEVERAL REPRESENTATIVE GROUPS OF DGS OPTIMAL CONFIGURATION

| Solutions | Optimal DG Capacity (bus locations to install DG) | | | | Objective Functions | | |
|---|---|---|---|---|---|---|---|
| | WT/kW (61) | PV/kW (50) | MT/kW (49) | MT/kW (64) | Line losses /kW | Investment benefit /pu | $VSF_{total}$ /pu |
| **solution A** | **1556** | **183** | **185** | **30** | **83.10** | **1.23** | **67.02** |
| solution B | 1515 | 122.6 | 130 | 73 | 83.56 | 1.20 | 67.01 |
| solution C | 1591 | 167 | 114.3 | 32 | 83.76 | 1.23 | 67.00 |
| solution D | 1562.7 | 124.7 | 126.1 | 13.4 | 83.85 | 1.22 | 66.99 |
| solution E | 1569 | 187 | 156.6 | 60 | 84.95 | 1.23 | 66.97 |

From Fig. 4, we can see that the distribution of the Pareto optimal solution is uniform and completed, which verifies the effectiveness of the proposed algorithm. The bus numbers of the corresponding installation locations are 61, 50, 49, and 64. From Table IV, we can find that, for the objective functions, the line losses, investment benefit, and the system voltage stability factor cannot achieve optimal results at the same time, and appropriate solutions can be chosen in practical applications according to the actual requirements.

In order to validate the effectiveness of the presented method, a comparison has been made between the 'best' comprise solution obtained by the proposed method and those obtained by other current state-of-the-art approaches, such as multi-objective harmony search (MOHS) [44], MOPSO [45] and NSGA-II [46]. The performance comparison of these schemes is shown in Table V.

TABLE V.   RESULT COMPARISON OF PG&E 69-BUS SYSTEM

| Objective functions | Base case without DGs | NSGA-II | MOPSO | MOHS | Proposed approach |
|---|---|---|---|---|---|
| $C_p$ (p.u.) | / | 1.01 | 1.13 | 1.21 | **1.23** |
| $VSF_{total}$ (p.u.) | 64.68 | 66.57 | 66.44 | 66.86 | **67.02** |
| Line losses (kW) | 224.94 | 96.73 | 89.56 | 88.16 | **83.10** |

From Table V it can be observed that: (1) compared with the base case without DGs, all the DG planning methodologies achieved good results. Especially, when using the proposal, line losses have been reduced by 63.06%, and voltage stability factor has been increased by 3.62%; (2) more importantly, that the proposed method is superior to all other approaches, embodying that its test results are the best ones in each objective function to be optimized. For example, compared with NSGA-II, the line losses of the proposal has been reduced by 14.09%, and DG investment benefits and voltage stability factor have been increased by 21.78% and 0.68%, respectively.

The above results prove that the proposed two-stage optimization method is more effective than other planning approaches. The main reason for this is that through an alternating iteration between the upper- and lower- layers, the optimization process can be guided towards the global Pareto-optimal front more effectively and efficiently in the search space, besides the strong global searching ability of the MOALO. Therefore, the conclusion can be drawn that the proposal is an effective way for solving optimal DG placement.

*D. Impact analysis of energy storage integration*

For ease of analysis, the optimal DG capacity of the best compromise solution is employed as the pre-designed rated capacities of DGs in this section. In this case, the probability density distributions of the power outputs of DGs, including WT and PV, are shown in Fig. 5.

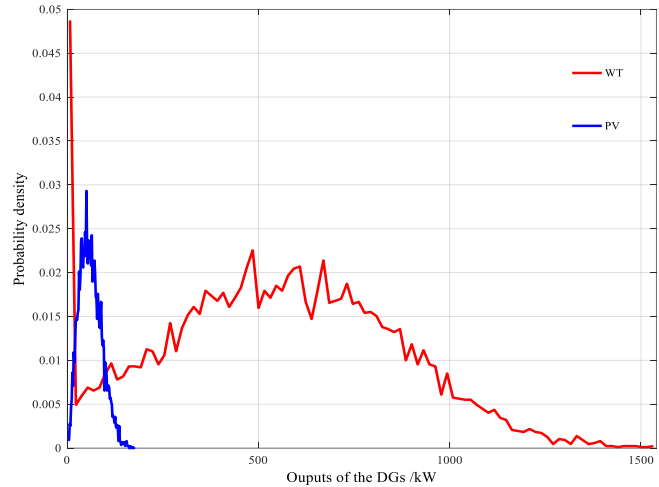

**Fig. 5.** Probability density distributions of the power outputs of DGs

Fig. 5 shows that the power outputs of renewable DG units, whether WT or PV, are unsatisfactory, since their outputs are much less than their expected rated capacities for most of the time. To solve this problem, ESDs are respectively installed at the buses 50 and 61, where the above-mentioned DGs are installed.

Once given a pre-assigned confidence level, the reserve energy storage output can be solved according to (14). The reserve energy storage outputs for different confidence levels are listed in Table VI.

TABLE VI.   RESERVE ENERGY STORAGE OUTPUT FOR DIFFERENT CONFIDENCE LEVELS

| Confidence | Reserve Output of Energy Storage / kW |
|---|---|

| Level | bus 61 | bus 50 |
|---|---|---|
| 95% | 1490 | 130 |
| 90% | 1190 | 110 |
| 85% | 990 | 90 |
| 80% | 820 | 80 |
| 75% | 680 | 70 |
| 70% | 560 | 60 |
| 65% | 450 | 50 |
| 60% | 360 | 40 |

Table VI demonstrates that the confidence level can be significantly improved with more energy storage capacity installed. But, on the other hand, this inevitably leads to that the economics decline to a certain extent due to the increasement of energy storage costs. For this reason, it is critical to choose an appropriate confidence level for achieving a reasonable tradeoff between the reliability and economics.

In order to verify the effects of the integration of energy storage, the cumulative probability distributions of the power outputs of DGs with/without ESDs are shown in Fig. 6.

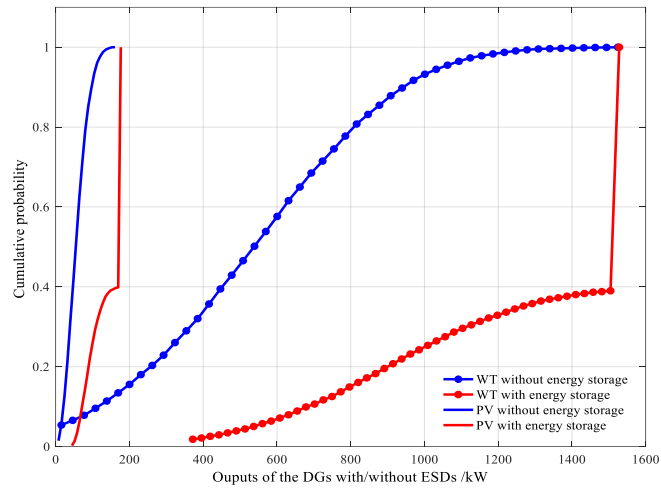

**Fig. 6.** Probability density distributions of the power outputs of DGs

It can be clearly seen from Fig. 6 that the power outputs of DGs are remarkably enhanced due to the integration of ESDs. Take, for example, power output of WT. Without ESDs, the power output of the bus connected to WT can only reach 800kW or below with a probability of 80%, and, especially, the power output probability is at the level of 1% or blew when reaching the rated capacity. After installing ESDs, the minimum power output of the bus connected to WT has been increased by 360 kW, and the probability that the power output reaches the rated capacity has been markedly improved to about 60%. The power output of PV is also similar to that of WT.

We calculate probabilistic power flow for the case with or without energy storage. The Monte Carlo method is then used to calculate the cumulative probability distribution curves of annual (8700 hours) power loss and voltage magnitude.

At the 60% confidence level, the cumulative probability distribution curves of the active and reactive power losses are shown in Figs. 7 and 8.

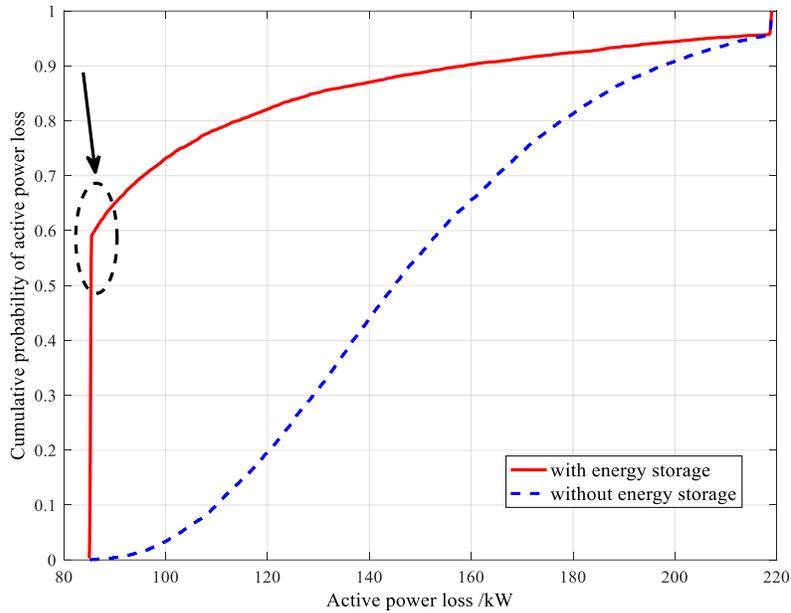

**Fig. 7.** Cumulative probability of active power loss

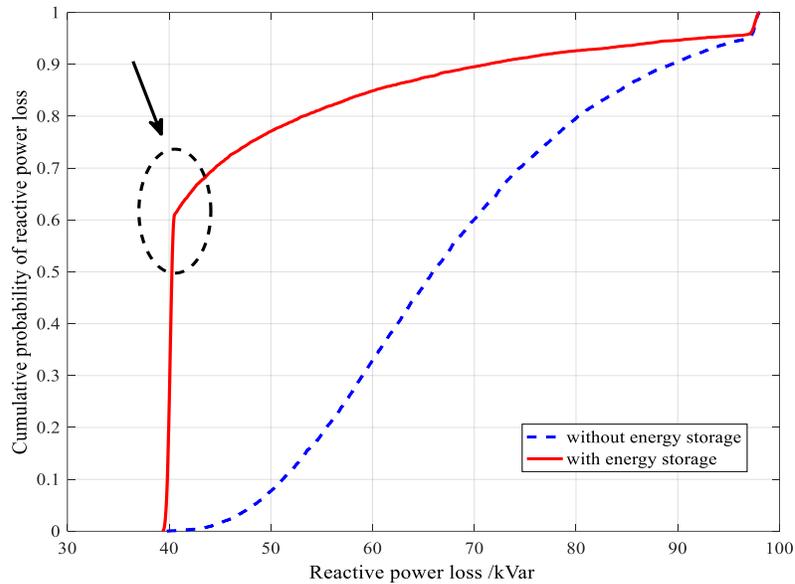

**Fig. 8.** Cumulative probability of reactive power loss

As is shown in Figs. 7 and 8, after the addition of energy storage, the cumulative probability distribution of the system power loss increases rapidly to 60% at the very beginning, and then begins to rise gradually. It indicates that the installed ESDs can reduce the system power losses to the level corresponding to solution A with a high probability. And thereby, the results prove that energy storage plays an important role in improving the overall economics of the system.

Since the minimum voltage in the system is at the bus 27, we perform a probabilistic simulation of the voltage magnitude of the bus 27. The cumulative probability distribution is shown in Fig. 9.

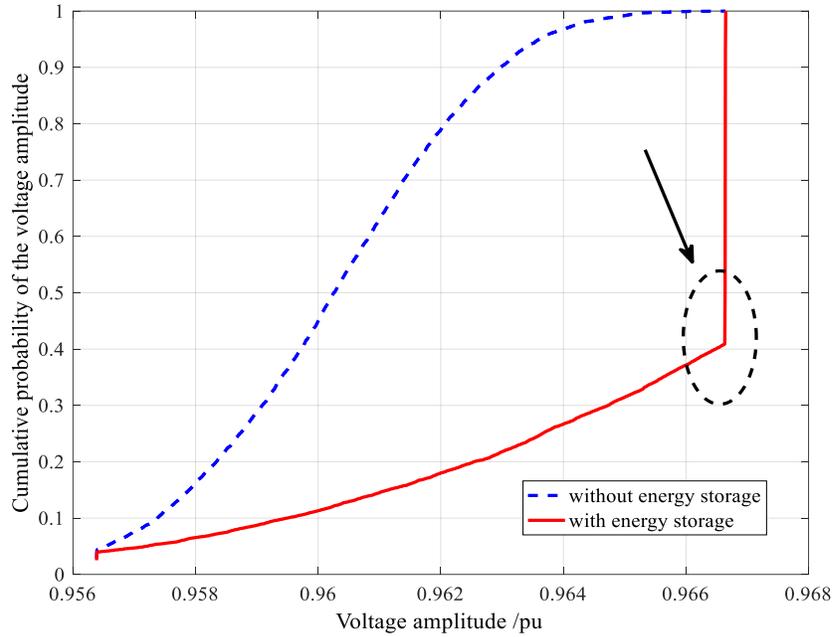

**Fig. 9.** Cumulative probability of the voltage amplitude at the bus 27

From Fig. 9, it is clearly seen that the voltage magnitude at the bus 27 can be stabilized between 0.967 pu~0.968 pu with a high probability (at least 60%) due to the integration of energy storage. By the same token, the voltage magnitudes at all other buses, except the swing bus, have similar results. Therefore, the conclusion can be drawn on the basis of the evidences that integration of energy storage is an effective and feasible way to improve the power output performances of DGs, which makes DGs operate more closely to their pre-designed rated capacities at the planning stage.

## V. CONCLUSION

Due to the intermittent nature of renewable DGs, traditional planning methods are increasingly less able to meet the needs of modern active distribution networks. Furthermore, the actual power outputs of DGs are hardly achieving their pre-designed rated capacities in the planning stage. To handle this problem, a new two-stage optimal DG planning method with the consideration of the integration of energy storages is presented in this paper. Test results on the PG&E 69-bus distribution system are displayed to verify the effectiveness of our new approach. As a result of this study, the following general conclusions can be drawn:

1) The proposed two-stage optimization method is superior to current state-of-the-art DG planning approaches, such as NSGA-II, MOPSO and MOHS, embodying that our method yields better results for all of the objectives such as higher investment benefits, better voltage stability and lower line losses.

2) The impact analysis of energy storage integration demonstrates that energy storage is an effective and feasible way to improve the power output performances of renewable DGs, which makes the DGs operate at their pre-designed rated capacities at the planning stage with the probability of at least 60%.

As a future work, the proposed methodology may find potential applications in an integrated planning of smart ADNs facing severe uncertainties resulting from both the generation side and the load side. In addition, our next step work will focus on considering to handle uncertainties of ADNs through the use of hybrid energy storage systems during the optimization process, since it allows the obtained planning schemes to be more committed with reality.

## REFERENCES

[1] Wang C, Song G, Li P, Ji H, Zhao J, Wu J. Optimal siting and sizing of soft open points in active electrical distribution networks. Appl Energy 2017; 189: 301-309.
[2] Jia H, Qi W, Liu Z, Wang B, Zeng Y, Xu T. Hierarchical risk assessment of transmission system considering the influence of active distribution network. IEEE Trans Power Syst 2015; 30: 1084-1093.
[3] Jin X, Mu Y, Jia H, Wu J, Jiang T, Yu X. Dynamic economic dispatch of a hybrid energy microgrid considering building based virtual energy storage system. Appl Energy 2017; 194: 386-398.
[4] Zhou Y, Wang C, Wu J, Wang J, Cheng M, Li G. Optimal scheduling of aggregated thermostatically controlled loads with renewable generation in the intraday electricity market. Appl Energy 2017; 188: 456-465.
[5] Soroudi A. Possibilistic-scenario model for DGs impact assessment on distribution networks in an uncertain environment. IEEE Trans Power Syst 2012; 27: 1283-1293.
[6] Keane A, Ochoa LF, Borges CLT, Ault GW, Alarcon-Rodriguez AD, Currie RAF, et al. State-of-the-art techniques and challenges ahead for distributed generation planning and optimization. IEEE Trans Power Syst 2013; 28(2): 1493-1502.
[7] Georgilakis PS, Hatziargyriou ND. Optimal distributed generation placement in power distribution networks: models, methods, and future research. IEEE Trans Power Syst 2013; 28(3): 3420-3428.
[8] Kazemi MA, Sedighizadeh M, Mirzaei MJ, Homaee O. Optimal siting and sizing of distribution system operator owned EV parking lots. Appl Energy 2016; 179: 1176-1184.
[9] Mokgonyana L, Zhang J, Li H, Hu Y. Optimal location and capacity planning for distributed generation with independent power production and self-generation. Appl Energy 2017; 188: 140-150.
[10] Wang C, Nehrir MH. Analytical approaches for optimal placement of distributed generation sources in power systems. IEEE Trans Power Syst 2004; 19(4): 2068-2076.
[11] Acharya N, Mahat P, Mithulananthan N. An analytical approach for DG allocation in primary distribution network. Int J Electr Power Energy Syst 2006; 28: 669-678.
[12] Hung DQ, Mithulananthan N. Multiple distributed generator placement in primary distribution networks for loss reduction. IEEE Trans Ind Electron 2013; 60(4): 1700-1708.
[13] Viral R, Khatod DK. An analytical approach for sizing and siting of DGs in balanced radial distribution networks for loss minimization. Int J Electr Power Energy Syst 2015; 67: 191-201.
[14] Atwa YM, El-Saadany EF, Salama MMA, Seethapathy R. Optimal renewable resources mix for distribution system energy loss minimization. IEEE Trans Power Syst 2010; 25(1): 360-370.
[15] Abri RSA, El-Saadany EF, Atwa YM. Optimal placement and sizing method to improve the voltage stability margin in a distribution system using distributed generation. IEEE Trans Power Syst 2013; 28(1): 326-334.
[16] Shaaban MF, Atwa YM, El-Saadany EF. DG allocation for benefit maximization in distribution networks. IEEE Trans Power Syst 2013; 28(2): 639-649.
[17] Lee SH, Park JW. Selection of optimal location and size of multiple distributed generations by using Kalman filter algorithm. IEEE Trans Power Syst 2009; 24(3): 1393-1400.
[18] Martins VF, Borges CLT. Active distribution network integrated planning incorporating distributed generation and load response uncertainties. IEEE Trans Power Syst 2011; 26: 2164-2172.
[19] Aman MM, Jasmon GB, Bakar AHA, Mokhlis H. A new approach for optimum simultaneous multi-DG distributed generation units placement and sizing based on maximization of system loadability using HPSO (hybrid particle swarm optimization) algorithm. Energy 2014; 66: 202-215.
[20] Kanwar N, Gupta N, Niazi KR. Swarnkar A, Bansal RC. Simultaneous allocation of distributed energy resource using improved particle swarm optimization. Appl Energy 2017; 185: 1684-1693.
[21] Mostafa N, Rachid C, Mario P. Optimal allocation of dispersed energy storage systems in active distribution networks for energy balance and grid support. IEEE Trans Power Syst 2014; 29: 2300-2310.
[22] Sedghi M, Ahmadian A, Aliakbar-Golkar M. Optimal storage planning in active distribution network considering uncertainty of wind power distributed generation. IEEE Trans Power Syst 2016; 31: 304-316.
[23] Santos SF, Fitiwi DZ, Cruz MRM, Cabrita CM, Catalão JP. Impacts of optimal energy storage deployment and network reconfiguration on renewable integration level in distribution systems. Appl Energy 2017; 185: 44-55.

[24]    Mahani K, Farzan F, Jafari MA. Network-aware approach for energy storage planning and control in the network with high penetration of renewables. Appl Energy 2017; 195: 974-990.
[25]    Tang Y, Low SH. Optimal placement of energy storage in distribution networks. IEEE Trans Smart Grid; to be appeared, doi: 10.1109/TSG.2017.2711921.
[26]    Luo X, Wang J, Dooner M, Clarke J. Overview of current development in electrical energy storage technologies and the application potential in power system operation. Appl Energy 2015; 137: 511-536.
[27]    Dunn B, Kamath H, Tarascon JM. Electrical energy storage for the grid: a battery of choices. Science 2011; 334(6058), 928-935.
[28]    Cakici M, Kakarla RR, Alonso-Marroquin F. Advanced electrochemical energy storage supercapacitors based on the flexible carbon fiber fabric-coated with uniform coral-like $MnO_2$, structured electrodes. Chem Eng J 2017; 309: 151-158.
[29]    Ma Y, Chang H, Zhang M. Graphene-based materials for lithium-ion hybrid supercapacitors. Advanced Materials 2015; 27: 5296-5308.
[30]    Hassan M, Haque E, Reddy KR. Edge-enriched graphene quantum dots for enhanced photo-luminescence and super capacitance. Nanoscale 2014; 6: 11988-11994.
[31]    Hassan M, Reddy K R, Haque E. Hierarchical assembly of grapheme/polyaniline nanostructures to synthesize free-standing supercapacitor electrode. Compos Sci Technol 2014; 98: 1-8.
[32]    Bonaccorso F, Colombo L, Yu G, Stoller M, Tozzini V, Ferrari AC, et al. Graphene, related two-dimensional crystals, and hybrid systems for energy conversion and storage. Science 2015; 347(6217): 1246501.
[33]    Reddy KR, Hassan M, Gomes VG. Hybrid nanostructures based on titanium dioxide for enhanced photocatalysis. Applied Catalysis A General 2015; 489: 1-16.
[34]    Reddy KR, Sin BC, Ryu KS, Noh J, Lee Y. In situ self-organization of carbon black-polyaniline composites from nanospheres to nanorods: Synthesis, morphology, structure and electrical conductivity. Synthetic Metals 2009; 159(19): 1934-1939.
[35]    Reddy KR, Sin BC, Chi HY. A new one-step synthesis method for coating multi-walled carbon nanotubes with cuprous oxide nanoparticles. Scripta Materialia 2008; 58(11):1010-1013.
[36]    El-Khattam W, Bhattacharya K, Hegazy Y, Salama MMA. Optimal investment planning for distributed generation in a competitive electricity market. IEEE Trans Power Syst 2004; 19(3):1674–84.
[37]    Hung DQ, Mithulananthan N, Bansal RC. An optimal investment planning framework for multiple distributed generation units in industrial distribution systems. Appl Energy 2014; 124: 62-72.
[38]    Mirjalili S. The Ant Lion Optimizer. Adv Eng Softw 2015; 83: 80-98.
[39]    Ali ES, Elazim SMA, Abdelaziz AY. Ant lion optimization algorithm for renewable distributed generations. Energy 2016; 116: 445-458.
[40]    Dubey HM, Pandit M, Panigrahi BK. Ant lion optimization for short-term wind integrated hydrothermal power generation scheduling. Int J Electr Power Energy Syst 2016; 83: 158-174.
[41]    Mirjalili S, Jangir P, Saremi S. Multi-objective ant lion optimizer: a multi-objective optimization algorithm for solving engineering problems. Appl Intell 2017; 46: 79-95.
[42]    Zheng G, Jing Y, Huang H, Gao Y. Application of improved grey relational projection method to evaluate sustainable building envelope performance. Appl Energy 2010; 87: 710-720.
[43]    Ortega-Vazquez MA, Kirschen DS. Estimating the spinning reserve requirements in systems with significant wind power generation penetration. IEEE Trans Power Syst 2009; 24(1): 114-124.
[44]    Nekooei K, Farsangi MM, Nezamabadi-Pour H, Lee KY. An improved multi-objective harmony search for optimal placement of DGs in distribution systems. IEEE Trans Smart Grid 2013; 4(1): 557-567.
[45]    Zeinalzadeh A, Mohammadi Y, Moradi MH. Optimal multi objective placement and sizing of multiple DGs and shunt capacitor banks simultaneously considering load uncertainty via MOPSO approach. Int J Electr Power Energy Syst 2015; 67: 336-349.
[46]    Sheng W, Liu KY, Liu Y, Meng X, Li Y. Optimal placement and sizing of distributed generation via an improved nondominated sorting genetic algorithm II. IEEE Trans Power Del 2015; 30(2): 569-578.